\documentclass[11pt]{article}
\usepackage{amsfonts}
\usepackage{amsfonts}
\usepackage{amsfonts}
\usepackage{mathrsfs}
\usepackage{bm}
\usepackage{cite}
\usepackage[colorlinks,linkcolor=blue,citecolor=blue]{hyperref}
\usepackage{amssymb,amsmath} 
 \allowdisplaybreaks

\catcode`!=11
\let\!int\int \def\int{\displaystyle\!int}
\let\!lim\lim \def\lim{\displaystyle\!lim}
\let\!sum\sum \def\sum{\displaystyle\!sum}
\let\!sup\sup \def\sup{\displaystyle\!sup}
\let\!inf\inf \def\inf{\displaystyle\!inf}
\let\!cap\cap \def\cap{\displaystyle\!cap}
\let\!max\max \def\max{\displaystyle\!max}
\let\!min\min \def\min{\displaystyle\!min}
\let\!frac\frac \def\frac{\displaystyle\!frac}
\catcode`!=12

\let\oldsection\section
\renewcommand\section{\setcounter{equation}{0}\oldsection}

\allowdisplaybreaks
\def\pf{\it{Proof.}\rm\quad}

\def\N{\mathbb{N}}\def\Z{\mathbb{Z}}

\newtheorem{thm}{Theorem}[section]

\newtheorem{cor}[thm]{Corollary}

\begin{document}
\title {\bf Explicit Evaluations of Sums of Sequence Tails}
\author{
{Ce Xu\thanks{Corresponding author. Email: 15959259051@163.com (C. Xu)}\quad\  Xiaolan Zhou\thanks{Email: xlzhou@stu.xmu.edu.cn (X. Zhou)}}\\[1mm]
\small School of Mathematical Sciences, Xiamen University\\
\small Xiamen
361005, P.R. China}

\date{}
\maketitle \noindent{\bf Abstract }
In this paper, we use Abel's summation formula to evaluate several quadratic and cubic sums of the form:
\[{F_N}\left( {A,B;x} \right) := \sum\limits_{n = 1}^N {\left( {A - {A_n}} \right)\left( {B - {B_n}} \right){x^n}} ,\;x \in [ - 1,1]\]
and
\[F\left( {A,B,\zeta (r)} \right): = \sum\limits_{n = 1}^\infty  {\left( {A - {A_n}} \right)\left( {B - {B_n}} \right)\left( {\zeta \left( r \right) - {\zeta _n}\left( r \right)} \right)} ,\]
where the sequences $A_n,B_n$ are defined by the finite sums
 ${A_n} := \sum\limits_{k = 1}^n {{a_k}} ,\ {B_n} := \sum\limits_{k = 1}^n {{b_k}}\ ( {a_k},{b_k} =o(n^{-p}),{\mathop{\Re}\nolimits} \left( p \right) > 1 $)
 and $A = \mathop {\lim }\limits_{n \to \infty } {A_n},B = \mathop {\lim }\limits_{n \to \infty } {B_n},F\left( {A,B;x} \right) = \mathop {\lim }\limits_{n \to \infty } {F_n}\left( {A,B;x} \right)$. Namely, the sequences $A_n$ and $B_n$ are the partial sums of the convergent series $A$ and $B$, respectively.

We give an explicit formula of ${F_n}\left( {A,B;x} \right)$ by using the method of Abel's summation formula. Then we
use apply it to obtain a family of identities relating harmonic numbers to multiple zeta values. Furthermore, we also evaluate several other series involving multiple zeta star values. Some interesting (known or new) consequences and illustrative examples are considered.
\\[2mm]
\\[2mm]
\tableofcontents

\section{Introduction}
In recent years there has been significant interest in investigating sums of products of Riemann zeta tails and multiple zeta (star) values.  The subject of this paper is studying quadratic and cubic sums
\begin{align}\label{1.1}
F\left( {A,B;x} \right) := \sum\limits_{n = 1}^\infty  {\left( {A - {A_n}} \right)\left( {B - {B_n}} \right){x^n}},
\end{align}
\begin{align}\label{1.2}
F\left( {A,B,\zeta(r) } \right) := \sum\limits_{n = 1}^\infty  {\left( {A - {A_n}} \right)\left( {B - {B_n}} \right)\left( {\zeta \left( r \right) - {\zeta _n}\left( r \right)} \right)},\end{align}
are classical harmonic number and generalized harmonic numbers or partial sums of the series (see \cite{A2000})
\begin{align}\label{1.3}
\zeta \left( s \right):= \sum\limits_{j = 1}^\infty  {\frac{1}{{{j^s}}}} ,\;\Re(s)>1
\end{align}
defining the Riemann zeta function, respectively.
Here $\zeta_n(r)$ denotes the harmonic number (also called the partial sum of the Riemann zeta function), which is defined by
\begin{align}\label{1.4}
\zeta_n(r):=\sum\limits_{j=1}^n\frac {1}{j^r}\quad (r,n\in \N:=\{1,2,3\cdots\}),
\end{align}
with $H_n:=\zeta_n(1)$ is the classical harmonic number, and the empty sum $\zeta_0(r)$ is conventionally understood to be zero.
In recent papers \cite{OF2015,OF2016}, O. Furdui, A. S\^{i}nt\v{a}m\v{a}rian  and C. V\u{a}lean prove some results on sums of products of the ¡°tails¡± of the Riemann zeta function, i.e.,
\[\zeta \left( s \right) - {\zeta _n}\left( s \right) = \sum\limits_{j = n + 1}^\infty  {\frac{1}{{{j^s}}}}.\]
In below, we let $F\left( {A,B} \right): = F\left( {A,B;1} \right)$. In \cite{OF2015}, O. Furdui and A. S\^{i}nt\v{a}m\v{a}rian shown that quadratic series
\[F(\zeta(s),\zeta(s))=\sum\limits_{n = 1}^\infty  {{{\left( {\zeta \left( s \right) - {\zeta _n}\left( s \right)} \right)}^2}}, \ s\in \N\setminus\{1\}:=\{2,3,\ldots\}\]
can be expressed as a rational linear combination of zeta
function values, binomial coefficients and products of two zeta function values. In \cite{OF2016}, O. Furdui and C. V\u{a}lean proven that the series
\[F(\zeta(s),\zeta(s+1))=\sum\limits_{n = 1}^\infty  {\left( {\zeta \left( s \right) - {\zeta _n}\left( s \right)} \right)\left( {\zeta \left( {s + 1} \right) - {\zeta _n}\left( {s + 1} \right)} \right)}\quad (s\in \N\setminus\{1\}) \]
can be expressed in terms of Riemann zeta values and a special integral involving a polylogarithms. Here the polylogarithm function is defined as follows
\begin{align}\label{1.5}
{\rm Li}{_p}\left( x \right) := \sum\limits_{n = 1}^\infty  {\frac{{{x^n}}}{{{n^p}}}}, \Re(p)>1,\ \left| x \right| \le 1 .
\end{align}
Furthermore, in \cite{MEH2016}, Hoffman shown a simple general result for the quadratic sum
\[F(\zeta(p),\zeta(q))=\sum\limits_{n = 1}^\infty  {\left( {\zeta \left( p \right) - {\zeta _n}\left( p \right)} \right)\left( {\zeta \left( q \right) - {\zeta _n}\left( q \right)} \right)} ,\;p,q \in \N\setminus\{1\}\]
in terms of the multiple zeta values $\zeta \left( {{s_1},{s_2}, \ldots ,{s_m}} \right)$ (for short mzvs) defined by
\begin{align}\label{1.6}
\zeta \left( {{s_1},{s_2}, \ldots ,{s_m}} \right) :=  \sum\limits_{{k_1} >  \cdots  > {k_m} \ge 1} {\frac{1}{{k_1^{{s_1}} \cdots k_m^{{s_m}}}}},
\end{align}
where for convergence ${s_1} + {s_2} +  \cdots  + {s_j} > j$ for $j=1,2,\ldots,m$. Similarly, the the multiple zeta star values $\zeta^\star \left( {{s_1},{s_2}, \ldots ,{s_m}} \right)$ (for short mzsvs) defined by
\begin{align}
\zeta \left( {{s_1},{s_2}, \ldots ,{s_m}} \right) :=  \sum\limits_{{k_1} \geq  \cdots \geq {k_m} \ge 1} {\frac{1}{{k_1^{{s_1}} \cdots k_m^{{s_m}}}}},
\end{align}
where for convergence ${s_1} + {s_2} +  \cdots  + {s_j} > j$ for $j=1,2,\ldots,m$.
Moreover, Hoffman also provided an explicit evaluation of sums of products of three or more Riemann zeta tails in a closed form in terms of multiple zeta values (see Theorem 4.3 in the reference \cite{MEH2016}). The multiple zeta value $\zeta \left( {{s_1},{s_2}, \ldots ,{s_m}} \right)$ is said to have depth $m$ and weight $s_1+\ldots+s_m$. Multiple zeta values of general depth were introduced independently by Hoffman \cite{MEH1992} and D. Zagier \cite{DZ1994}, but for depth two they were already studied by Euler. In \cite{BG1996}, Borwein and Girgensohn gave some evaluation of multiple zeta values with depth three, and they proved that all $\zeta \left( {q,p,r} \right)$ with $r+p+q$ is even or less than or equal to 10 or $r+p+q=12$ were reducible to zeta values and double zeta values. Further, Eie and Wei \cite{EW2012} proved that for positive integer $p>1$, all mzvs $\zeta \left( {p,p,1,1} \right)$ can be expressed in terms of mzvs of depth $\leq 3$, and gave explicit formulas. There are also a lot of contributions on multiple zeta values in the last two decades, see \cite{BBBL1997,BBBL2001,MEH1992,MEH1997,Xu2017,DZ1994,DZ2012} and references therein.

In general, the alternating multiple zeta values and multiple zeta star values are defined by
\begin{align}\label{1.7}
\zeta \left( \bf s \right)\equiv\zeta \left( {{s_1}, \ldots ,{s_m}} \right): = \sum\limits_{{k_1} >  \cdots  > {k_m} \geq 1} {\prod\limits_{j = 1}^m {n_j^{ - \left| {{s_j}} \right|}}{\rm sgn}(s_j)^{k_j}},
\end{align}
\begin{align}\label{1.8}
\zeta^\star \left( \bf s \right)\equiv{\zeta ^ \star }\left( {{s_1}, \ldots ,{s_m}} \right): = \sum\limits_{{k_1} \ge  \cdots  \ge {k_m} \ge 1} {\prod\limits_{j = 1}^m {n_j^{ - \left| {{s_j}} \right|}{\rm sgn}(s_j)^{k_j}} },
\end{align}
where ${\bf s}:=\{s_1,\ldots,s_m\}\in (\Z^*)^m\ (\Z^*:=\{\pm1,\pm2,\ldots\})$ with $s_1\neq 1$ , and
\[{\mathop{\rm sgn}} \left( {{s_j}} \right): = \left\{ {\begin{array}{*{20}{c}}
   {1,} & {{s_j} > 0,}  \\
   { - 1,} & {{s_j} < 0.}  \\
\end{array}} \right.\]
Throughout the paper we will use the notation $\bar p$ to denote a negative entry $s_j=-p$\ ($p$ is a positive integer). For example,
 \[\zeta \left( {{{\bar s}_1},{s_2}} \right) = \zeta \left( { - {s_1},{s_2}} \right),\quad {\zeta ^ \star }\left( {{{\bar s}_1},{s_2},{{\bar s}_3}} \right) = {\zeta ^ \star }\left( { - {s_1},{s_2}, - {s_3}} \right).\]
We call $l({\bf s}):=m$ and $\left| {\bf s} \right|: = \sum\limits_{j = 1}^k {\left| {{s_j}} \right|} $ the depth and the weight of (\ref{1.7}) and (\ref{1.8}), respectively. For convenience we let $\zeta \left( \emptyset  \right) = {\zeta ^ \star }\left( \emptyset  \right) = 1$ and ${\left\{ {{s_1}, \ldots ,{s_j}} \right\}_d}$ the set formed by repeating the composition $\left( {{s_1}, \ldots ,{s_j}} \right)$ $d$ times.

Furthermore, for convenience, we define the generalized multiple harmonic and star sums by
\begin{align*}
&\zeta_n \left( \bf s \right)\equiv\zeta_n \left( {{s_1}, \ldots ,{s_m}} \right): = \sum\limits_{n\geq {k_1} >  \cdots  > {k_m} \geq 1} {\prod\limits_{j = 1}^m {n_j^{ - \left| {{s_j}} \right|}}{\rm sgn}(s_j)^{k_j}},\\
&\zeta_n^\star \left( \bf s \right)\equiv{\zeta_n ^ \star }\left( {{s_1}, \ldots ,{s_m}} \right): = \sum\limits_{n\geq {k_1} \ge  \cdots  \ge {k_m} \ge 1} {\prod\limits_{j = 1}^m {n_j^{ - \left| {{s_j}} \right|}{\rm sgn}(s_j)^{k_j}} },
\end{align*}
which are also called the partial sums of the multiple zeta and zeta star values, respectively. Here ${\bf s}:=\{s_1,\ldots,s_m\}\in (\Z^*)^m$

In this paper, we discuss the analytic representations of $F(A,B;x)$ and $F(A,B,\zeta(r))$. The results which we present here can be seen as an extension of Hoffman's, Furdui's and V\u{a}lean's work. Many results of this paper can be expressed as a rational linear combination of products of zeta values and polylogarithms.
Using certain integral representations of series, we can prove that the cubic series
\[F\left( {\zeta \left( p \right),\zeta \left( q \right),\zeta \left( r \right)} \right) = \sum\limits_{n = 1}^\infty  {\left( {\zeta \left( p \right) - {\zeta _n}\left( p \right)} \right)\left( {\zeta \left( q \right) - {\zeta _n}\left( q \right)} \right)\left( {\zeta \left( r \right) - {\zeta _n}\left( r \right)} \right)}\quad (p,q,r \in \N\setminus\{1\})\]
can be expressed in terms of multiple zeta values or multiple zeta star values with depth $\leq 3$.
Furthermore, we also give some evaluation of sums of products of alternating Riemann zeta Tails. For example, we show that
\begin{align}\label{1.9}
&\sum\limits_{n = 1}^\infty  {\left( {\zeta \left( m \right) - {\zeta _n}\left( m \right)} \right)\left( {\bar \zeta \left( p \right) - {{\bar \zeta }_n}\left( p \right)} \right)}\nonumber \\&  =  - {\zeta ^ \star }\left( {\bar p,m - 1} \right) - {\zeta ^ \star }\left( {m,\overline {p - 1}} \right) - \zeta \left( m \right)\bar \zeta \left( p \right) - \bar \zeta \left( {m + p - 1} \right),
\end{align}
where ${\bar \zeta \left( s \right)}$ and $ {{{\bar \zeta }_n}\left( s \right)}$ denote the alternating Riemann zeta function and its partial sum (or called the alternating harmonic number), respectively, which are defined by
\[\bar \zeta \left( s \right) := \sum\limits_{j= 1}^\infty  {\frac{{{{\left( { - 1} \right)}^{j - 1}}}}{{{j^s}}}} ,\ {{\bar \zeta }_n}\left( s \right) := \sum\limits_{j= 1}^n {\frac{{{{\left( { - 1} \right)}^{j - 1}}}}{{{j^s}}}} ,\quad n \in \N, \Re(s)\geq1. \]
From above definitions, we know that ${\bar \zeta}\left(s\right)=-\zeta\left(\bar s\right)$. Hence, the formula (\ref{1.9}) can be rewritten as
\begin{align}\label{1.10}
F\left( {\zeta \left( m \right),\bar \zeta \left( p \right)} \right)= \zeta \left( m \right)\zeta \left( {\bar p} \right) + \zeta \left( \overline {m + p - 1} \right) - {\zeta ^ \star }\left( {\bar p,m - 1} \right) - {\zeta ^ \star }\left( {m,\overline {p - 1}} \right),
\end{align}
Similarly, by using the above notations, we have the relations
\begin{align*}
 &{\zeta ^ \star }\left( {m,p} \right) = \sum\limits_{n = 1}^\infty  {\frac{{{\zeta _n}\left( p \right)}}{{{n^m}}}} ,{\zeta ^ \star }\left( {\bar m,\bar p} \right) = \sum\limits_{n = 1}^\infty  {\frac{{{{\bar \zeta }_n}\left( p \right)}}{{{n^m}}}{{\left( { - 1} \right)}^{n - 1}}} , \\
 &{\zeta ^ \star }\left( {\bar m,p} \right) =  - \sum\limits_{n = 1}^\infty  {\frac{{{\zeta _n}\left( p \right)}}{{{n^m}}}{{\left( { - 1} \right)}^{n - 1}}} ,{\zeta ^ \star }\left( {m,\bar p} \right) =  - \sum\limits_{n = 1}^\infty  {\frac{{{{\bar \zeta }_n}\left( p \right)}}{{{n^m}}}} .
\end{align*}
The evaluation of linear sums (or double sums)
 \[\sum\limits_{n = 1}^\infty  {\frac{{{\zeta _n}\left( p \right)}}{{{n^m}}}} ,\sum\limits_{n = 1}^\infty  {\frac{{{{\bar \zeta }_n}\left( p \right)}}{{{n^m}}}{{\left( { - 1} \right)}^{n - 1}}} ,\sum\limits_{n = 1}^\infty  {\frac{{{\zeta _n}\left(p \right)}}{{{n^m}}}{{\left( { - 1} \right)}^{n - 1}}} {\rm and} \sum\limits_{n = 1}^\infty  {\frac{{{{\bar \zeta }_n}\left( p \right)}}{{{n^m}}}}\]
 in terms of values of Riemann zeta function (or alternating Riemann zeta function) and polylogarithm
function at positive integers is known when  $(p,m) = (1,3), (2,2)$ or $p+m$ is odd (see \cite{BBG1994,BBG1995,FS1998,S2015,X2016,CX2016}). For example, we have
\begin{align*}
&\sum\limits_{n = 1}^\infty  {\frac{{{{\bar \zeta}_n}\left( 1 \right)}}{{{n^3}}}}  = \frac{7}{4}\zeta \left( 3 \right)\ln 2 - \frac{5}{{16}}\zeta (4),\\
&\sum\limits_{n = 1}^\infty  {\frac{{{H_n}{}}}{{{n^3}}}}{\left( { - 1} \right)}^{n - 1}  =  - 2{\rm Li}{_4}\left( {\frac{1}{2}} \right) + \frac{{11{}}}{{4}} \zeta(4) + \frac{{{1}}}{{2}}\zeta(2){\ln ^2}2 - \frac{1}{{12}}{\ln ^4}2 - \frac{7}{4}\zeta \left( 3 \right)\ln 2,\\
&\sum\limits_{n = 1}^\infty  {\frac{{{{\bar \zeta}_n}\left( 1 \right)}}{{{n^3}}}{{\left( { - 1} \right)}^{n - 1}} = \frac{3}{2}\zeta \left( 4 \right)}  + \frac{1}{2}\zeta \left( 2 \right){\ln ^2}2 - \frac{1}{{12}}{\ln ^4}2 - 2{\rm Li}{_4}\left( {\frac{1}{2}} \right),\\
&\sum\limits_{n = 1}^\infty  {\frac{{{\zeta _n}\left( 2 \right)}}{{{n^2}}}{{\left( { - 1} \right)}^{n - 1}}}  =  - \frac{{51}}{{16}}\zeta (4) + 4{\rm{L}}{{\rm{i}}_4}\left( {\frac{1}{2}} \right) + \frac{7}{2}\zeta (3)\ln 2 - \zeta (2){\ln ^2}(2) + \frac{1}{6}{\ln ^4}2.
\end{align*}
Furthermore, we consider the nested sum
\[\sum\limits_{n = 1}^\infty  {\frac{{{y^n}}}{{{n^m}}}\left( {\sum\limits_{k = 1}^n {\frac{{{x^k}}}{{{k^p}}}} } \right)} ,\;x,y \in \left[ { - 1,1} \right),\;m,p \in \N.\]
By taking the sum over complementary pairs of summation indices, we obtain a simple reflection formula
\[\sum\limits_{n = 1}^\infty  {\frac{{{y^n}}}{{{n^m}}}\left( {\sum\limits_{k = 1}^n {\frac{{{x^k}}}{{{k^p}}}} } \right)}  + \sum\limits_{n = 1}^\infty  {\frac{{{x^n}}}{{{n^p}}}\left( {\sum\limits_{k = 1}^n {\frac{{{y^k}}}{{{k^m}}}} } \right)}  = {\rm Li}{_p}\left( x \right){\rm Li}{_m}\left( y \right) + {\rm Li}{_{p + m}}\left( {xy} \right).\]
Therefore, we can deduce that
\begin{align*}
&\sum\limits_{n = 1}^\infty  {\frac{{{{\bar \zeta}_n}\left( 2 \right)}}{{{n^2}}}{{\left( { - 1} \right)}^{n - 1}} = } \frac{{13}}{{16}}\zeta \left( 4 \right),\\
&\sum\limits_{n = 1}^\infty  {\frac{{{\zeta _n}\left( 3 \right)}}{n}{{\left( { - 1} \right)}^{n - 1}}}  = \frac{{19}}{{16}}\zeta (4) - \frac{3}{4}\zeta \left( 3 \right)\ln 2,\\
&\sum\limits_{n = 1}^\infty  {\frac{{{{\bar \zeta}_n}\left( 2 \right)}}{{{n^2}}} = \frac{{85}}{{16}}\zeta \left( 4 \right)}  - 4{\rm{L}}{{\rm{i}}_4}\left( {\frac{1}{2}} \right) + \zeta \left( 2 \right){\ln ^2}2 - \frac{1}{6}{\ln ^4}2 - \frac{7}{2}\zeta \left( 3 \right)\ln 2,\\
&\sum\limits_{n = 1}^\infty  {\frac{{{{\bar \zeta}_n}\left( 3 \right)}}{n}{{\left( { - 1} \right)}^{n - 1}} = } 2{\rm{L}}{{\rm{i}}_4}\left( {\frac{1}{2}} \right) + \frac{3}{4}\zeta \left( 3 \right)\ln 2 + \frac{1}{{12}}{\ln ^4}2 - \frac{1}{2}\zeta \left( 4 \right) - \frac{1}{2}\zeta \left( 2 \right){\ln ^2}2.
\end{align*}
Some results for sums of harmonic numbers or alternating harmonic numbers may be seen in the works of \cite{BBG1994,BBG1995,FS1998,M2014,S2015,S2011,S2010,X2016,CX2016} and references therein.

The main results of this paper are the following theorems.
\begin{thm}\label{thm1.1} If ${a_k},{b_k} =o(n^{-p}), \Re(p)>1$ and $x\in [-1,1]$, then
\begin{align}\label{1.11}
&\left( {1 - x} \right){F_n}\left( {A,B;x} \right) = \left( {1 - x} \right)\sum\limits_{k = 1}^n {\left( {A - {A_k}} \right)\left( {B - {B_k}} \right){x^k}}\nonumber \\
& = ABx - \left( {A - {A_n}} \right)\left( {B - {B_n}} \right){x^{n + 1}} - \left( {\sum\limits_{k = 1}^n {{a_k}{x^k}} } \right)\left( {B - {B_n}} \right) - \left( {\sum\limits_{k = 1}^n {{b_k}{x^k}} } \right)\left( {A - {A_n}} \right)\nonumber\\
&\quad - \sum\limits_{k = 1}^n {\left( {\sum\limits_{i = 1}^k {{a_i}{x^i}} } \right){b_k}}  - \sum\limits_{k = 1}^n {\left( {\sum\limits_{i = 1}^k {{b_i}{x^i}} } \right){a_k}}  + \sum\limits_{k = 1}^n {{a_k}{b_k}{x^k}} .
\end{align}
\end{thm}
\begin{thm}\label{thm1.2}If ${a_k},{b_k} =o(n^{-p}), \Re(p)>1$ and $r\in \N\setminus\{1\}$, then
\begin{align}\label{1.12}
F\left( {A,B,\zeta(r) } \right) =& \sum\limits_{k = 1}^\infty  {\left( {A - {A_k}} \right)\left( {B - {B_k}} \right)\left( {\zeta \left( r \right) - {\zeta _k}\left( r \right)} \right)}\nonumber \\
=&\zeta \left( r \right)\left\{ {\sum\limits_{n = 1}^\infty  {\left( {\sum\limits_{k = 1}^n {k{a_k}} } \right){b_n}}  + \sum\limits_{n = 1}^\infty  {\left( {\sum\limits_{k = 1}^n {k{b_k}} } \right){a_n}}  - \sum\limits_{n = 1}^\infty  {n{a_n}{b_n}}  - AB} \right\}\nonumber\\
&- \sum\limits_{n = 1}^\infty  {\left( {\sum\limits_{k = 1}^n {{a_k}\left( {k{\zeta _k}\left( r \right) - {\zeta _k}\left( {r - 1} \right)} \right)} } \right){b_n}}  - \sum\limits_{n = 1}^\infty  {\left( {\sum\limits_{k = 1}^n {{b_k}\left( {k{\zeta _k}\left( r \right) - {\zeta _k}\left( {r - 1} \right)} \right)} } \right){a_n}} \nonumber\\
& + \sum\limits_{n = 1}^\infty  {{a_n}{b_n}\left( {n{\zeta _n}\left( r \right) - {\zeta _n}\left( {r - 1} \right)} \right)} .
\end{align}
\end{thm}

\section{Proofs of Main Theorems}
In this section, we will use Abel's summation formula to prove Theorem \ref{thm1.1} and \ref{thm1.2}. The Abel's summation by parts formula states that if $(a_n)_{n\geq1}$ and $(b_n)_{n\geq1}$ are two sequences of real numbers and $A_n=\sum\limits_{k = 1}^n {{a_k}}$, then \[\sum\limits_{k = 1}^n {{a_k}{b_k}}  = {A_n}{b_{n + 1}} + \sum\limits_{k = 1}^n {{A_k}\left( {{b_k} - {b_{k + 1}}} \right)}.\]
 Firstly, we are ready to prove Theorem \ref{thm1.1}.

\subsection{Proof of Theorem 1.1}
 From the definition of function $F_n(A,B;x)$, we have
\begin{align}\label{2.1}
{F_n}\left( {A,B;x} \right) &= \sum\limits_{k = 1}^n {\left( {A - {A_k}} \right)\left( {B - {B_k}} \right){x^k}}\nonumber \\
& = \left( {A - {a_1}} \right)\left( {B - {b_1}} \right)x + \sum\limits_{k = 1}^{n - 1} {\left( {A - {A_{k + 1}}} \right)\left( {B - {B_{k + 1}}} \right){x^{k + 1}}}\nonumber \\
& = \left( {A - {a_1}} \right)\left( {B - {b_1}} \right)x + \sum\limits_{k = 1}^{n - 1} {\left( {A - {A_k} - {a_{k + 1}}} \right)\left( {B - {B_k} - {b_{k + 1}}} \right){x^{k + 1}}}\nonumber \\
& = \left( {A - {a_1}} \right)\left( {B - {b_1}} \right)x + \sum\limits_{k = 1}^{n - 1} {\left( {A - {A_k}} \right)\left( {B - {B_k}} \right){x^{k + 1}}}\nonumber \\
&\quad - \sum\limits_{k = 1}^{n - 1} {\left\{ {{a_{k + 1}}\left( {B - {B_k}} \right) + {b_{k + 1}}\left( {A - {A_k}} \right)} \right\}{x^{k + 1}}}  + \sum\limits_{k = 1}^{n - 1} {{a_{k + 1}}{b_{k + 1}}{x^{k + 1}}}\nonumber \\
& = \left( {A - {a_1}} \right)\left( {B - {b_1}} \right)x + \sum\limits_{k = 1}^n {\left( {A - {A_k}} \right)\left( {B - {B_k}} \right){x^{k + 1}}}  - \left( {A - {A_n}} \right)\left( {B - {B_n}} \right){x^{n + 1}}\nonumber\\
&\quad + \sum\limits_{k = 1}^n {{a_k}{b_k}{x^k}}  - {a_1}{b_1}x - \sum\limits_{k = 1}^n {\left\{ {{a_k}\left( {B - {B_k} + {b_k}} \right) + {b_k}\left( {A - {A_k} + {a_k}} \right)} \right\}{x^k}}  + {a_1}B + {b_1}A\nonumber\\
& = x{F_n}\left( {A,B;x} \right) + ABx - \left( {A - {A_n}} \right)\left( {B - {B_n}} \right){x^{n + 1}}\nonumber\\
&\quad - \sum\limits_{k = 1}^n {\left\{ {{a_k}\left( {B - {B_k}} \right) + {b_k}\left( {A - {A_k}} \right)} \right\}{x^k}}  - \sum\limits_{k = 1}^n {{a_k}{b_k}{x^k}}.
\end{align}
Then by using the Abel's summation formula, we can find that
\begin{align}\label{2.2}
\sum\limits_{k = 1}^n {\left\{ {{a_k}\left( {B - {B_k}} \right)} \right\}{x^k}}  = \left( {\sum\limits_{k = 1}^n {{a_k}{x^k}} } \right)\left( {B - {B_n}} \right) + \sum\limits_{k = 1}^n {\left( {\sum\limits_{i = 1}^k {{a_i}{x^i}} } \right){b_k}}  - \sum\limits_{k = 1}^n {{a_k}{b_k}{x^k}},
\end{align}
\begin{align}\label{2.3}
\sum\limits_{k = 1}^n {\left\{ {{b_k}\left( {A - {A_k}} \right)} \right\}{x^k}}  = \left( {\sum\limits_{k = 1}^n {{b_k}{x^k}} } \right)\left( {A - {A_n}} \right) + \sum\limits_{k = 1}^n {\left( {\sum\limits_{i = 1}^k {{b_i}{x^i}} } \right){a_k}}  - \sum\limits_{k = 1}^n {{a_k}{b_k}{x^k}}.
\end{align}
Hence, substituting formulas (\ref{2.2}) and (\ref{2.3}) into (\ref{2.1}) respectively, we deduce the (\ref{1.11}). Thus, the proof of Theorem \ref{thm1.1} is finished.

\subsection{Proof of Theorem 1.2}
Now we are ready to prove Theorem \ref{thm1.2}. Multiplying (\ref{1.11}) by $\frac{{{{\ln }^{r - 1}}x}}{{{{\left( {1 - x} \right)}^2}}}$ and integrating over the interval (0,1), and using the following elementary integral identities
\begin{align*}
&\int\limits_0^1 {\frac{{{x^k}{{\ln }^{r - 1}}x}}{{1 - x}}} dx = {\left( { - 1} \right)^{r - 1}}\left( {r - 1} \right)!\left( {\zeta \left( r \right) - {\zeta _k}\left( r \right)} \right),\;r \in \N\setminus\{1\},\\
&\int\limits_0^1 {\frac{{{x^k}{{\ln }^{r - 1}}x}}{{{{\left( {1 - x} \right)}^2}}}} dx = {\left( { - 1} \right)^{r - 1}}\left( {r - 1} \right)!\left( {\zeta \left( {r - 1} \right) - {\zeta _k}\left( {r - 1} \right) - k\zeta \left( r \right) + k{\zeta _k}\left( r \right)} \right),r \in \N\setminus\{1,2\},
\end{align*}
 by a direct calculation, we obtain
\begin{align}\label{2.4}
F\left( {A,B,\zeta \left( r \right)} \right) =& AB\left( {\zeta \left( {r - 1} \right) - \zeta \left( r \right)} \right)\nonumber\\
& + \sum\limits_{k = 1}^\infty  {{a_k}{b_k}\left( {\zeta \left( {r - 1} \right) - {\zeta _k}\left( {r - 1} \right) - k\zeta \left( r \right) + k{\zeta _k}\left( r \right)} \right)}\nonumber \\
&- \sum\limits_{n = 1}^\infty  {\left( {\sum\limits_{k = 1}^n {{a_k}\left( {\zeta \left( {r - 1} \right) - {\zeta _k}\left( {r - 1} \right) - k\zeta \left( r \right) + k{\zeta _k}\left( r \right)} \right)} } \right){b_n}}\nonumber \\
&- \sum\limits_{n = 1}^\infty  {\left( {\sum\limits_{k = 1}^n {{b_k}\left( {\zeta \left( {r - 1} \right) - {\zeta _k}\left( {r - 1} \right) - k\zeta \left( r \right) + k{\zeta _k}\left( r \right)} \right)} } \right){a_n}} .
\end{align}
Taking $x=1$ in (\ref{1.11}), we get
\begin{align}\label{2.5}
\sum\limits_{k = 1}^n {{A_k}{b_k}}  + \sum\limits_{k = 1}^n {{B_k}{a_k}}  = \sum\limits_{k = 1}^n {{a_k}{b_k}}  + {A_n}{B_n}.
\end{align}
Letting $n\rightarrow \infty$, then
\begin{align}\label{2.6}
\sum\limits_{k = 1}^\infty  {{A_k}{b_k}}  + \sum\limits_{k = 1}^\infty  {{B_k}{a_k}}  = \sum\limits_{k = 1}^\infty  {{a_k}{b_k}}  + AB.
\end{align}
Substituting (\ref{2.6}) into (\ref{2.4}), we know that when $r \in \N\setminus\{1,2\}$, the equality (\ref{1.12}) holds.

Similarly, multiplying (\ref{1.11}) by $\frac{{{{\ln }}x}}{{{{\left( {1 - x} \right)}^2}}}$ and integrating over the interval $(0,t)$, then letting $t\rightarrow 1$, we can prove that if $r=2$, the formula (\ref{1.2}) is also true. The proof of Theorem \ref{thm1.2} is finished.
\section{Some Results on Zeta Tails}
In this section, we will give some closed form of sums of sequences tails.

Firstly, multiplying $(1-x)^{-1}$ in both sides of (\ref{1.11}), then letting $x\rightarrow 1$ and using L'Hopital's rule, we deduce that
\begin{align}\label{3.1}
&\sum\limits_{k = 1}^n {\left( {A - {A_k}} \right)\left( {B - {B_k}} \right)}  = {F_n}\left( {A,B;1} \right) = \mathop {\lim }\limits_{x \to 1} \frac{{{(1-x)F_n}\left( {A,B;x} \right)}}{{1 - x}}\nonumber\\
 =& \sum\limits_{k = 1}^n {\left( {\sum\limits_{i = 1}^k {i{a_i}} } \right){b_k}}  + \sum\limits_{k = 1}^n {\left( {\sum\limits_{i = 1}^k {i{b_i}} } \right){a_k}}  + \left( {\sum\limits_{k = 1}^n {k{a_k}} } \right)\left( {B - {B_n}} \right) + \left( {\sum\limits_{k = 1}^n {k{b_k}} } \right)\left( {A - {A_n}} \right)\nonumber
 \\& + \left( {n + 1} \right)\left( {A - {A_n}} \right)\left( {B - {B_n}} \right) - \left( {\sum\limits_{k = 1}^n {k{a_k}{b_k}} } \right) - AB.
\end{align}
Letting $n\rightarrow \infty$ in (\ref{3.1}),  we arrive at the conclusion that
\begin{align}\label{3.2}
\sum\limits_{n = 1}^\infty  {\left( {A - {A_n}} \right)\left( {B - {B_n}} \right)}  = \sum\limits_{n = 1}^\infty  {\left( {\sum\limits_{k = 1}^n {k{a_k}} } \right){b_n}}  + \sum\limits_{n = 1}^\infty  {\left( {\sum\limits_{k = 1}^n {k{b_k}} } \right){a_n}}  - \sum\limits_{n = 1}^\infty  {n{a_n}{b_n}}  - AB.
\end{align}
Putting $\left( {{a_n},{b_n}} \right) = \left( {\frac{1}{{{n^m}}},\frac{1}{{{n^p}}}} \right),\left( {\frac{{{{\left( { - 1} \right)}^{n }}}}{{{n^m}}},\frac{{{{\left( { - 1} \right)}^{n }}}}{{{n^p}}}} \right),\left( {\frac{1}{{{n^m}}},\frac{{{{\left( { - 1} \right)}^{n }}}}{{{n^p}}}} \right)$ in (\ref{3.2}), and using the following relations
\begin{align*}
 &{\zeta ^ \star }\left( {m,p} \right) = \zeta \left( {m,p} \right) + \zeta \left( {m + p} \right), \\
 &{\zeta ^ \star }\left( {\bar m,\bar p} \right) = \zeta \left( {\bar m,\bar p} \right) + \zeta \left( {m + p} \right), \\
 &{\zeta ^ \star }\left( {\bar m,p} \right) = \zeta \left( {\bar m,p} \right) + \zeta \left( \overline{m + p} \right), \\
 &{\zeta ^ \star }\left( {m,\bar p} \right) = \zeta \left( {m,\bar p} \right) + {\zeta }\left( \overline {m+p} \right),
\end{align*}
we obtain the following identities
\begin{align}\label{3.3}
F\left( {\zeta \left( m \right),\zeta \left( p \right)} \right)&=\sum\limits_{n = 1}^\infty  {\left( {\zeta \left( m \right) - {\zeta _n}\left( m \right)} \right)\left( {\zeta \left( p \right) - {\zeta _n}\left( p \right)} \right)}\nonumber \\
&=\zeta \left( {p,m - 1} \right) + \zeta \left( {m,p - 1} \right) + \zeta \left( {m + p - 1} \right) - \zeta \left( m \right)\zeta \left( p \right),
\end{align}
\begin{align}\label{3.4}
F\left( {\zeta \left( {\bar m} \right),\zeta \left( {\bar p} \right)} \right) &= \sum\limits_{n = 1}^\infty  {\left( {\zeta \left( {\bar m} \right) - {\zeta _n}\left( {\bar m} \right)} \right)\left( {\zeta \left( {\bar p} \right) - {\zeta _n}\left( {\bar p} \right)} \right)}  \nonumber\\
& = \zeta \left( {\bar m,\overline {p - 1}} \right) + \zeta \left( {\bar p,\overline{m - 1}} \right) + \zeta \left( {m + p - 1} \right) - \zeta \left( {\bar m} \right)\zeta \left( {\bar p} \right),
\end{align}
\begin{align}\label{3.5}
 F\left( {\zeta \left( m \right),\zeta \left( {\bar p} \right)} \right) &= \sum\limits_{n = 1}^\infty  {\left( {\zeta \left( m \right) - {\zeta _n}\left( m \right)} \right)\left( {\zeta \left( {\bar p} \right) - {\zeta _n}\left( {\bar p} \right)} \right)}  \nonumber\\
 & = \zeta \left( {m,\overline{p - 1}} \right) + \zeta \left( {{\bar p},m - 1} \right) + \zeta \left( \overline{m + p - 1} \right) - \zeta \left( m \right)\zeta \left( {\bar p} \right) .
\end{align}
Similarly, in the case $\left( {{a_n},{b_n}} \right) = \left( {\frac{1}{{{n^m}}},\frac{1}{{{n^p}}}} \right)$, Theorem \ref{thm1.2} gives
\begin{align}\label{3.6}
F\left( {\zeta \left( m \right),\zeta \left( p \right),\zeta \left( r \right)} \right) = &\sum\limits_{n = 1}^\infty  {\left( {\zeta \left( m \right) - {\zeta _n}\left( m \right)} \right)\left( {\zeta \left( p \right) - {\zeta _n}\left( p \right)} \right)\left( {\zeta \left( r \right) - {\zeta _n}\left( r \right)} \right)}\nonumber \\
=&  \zeta \left( r \right)\left\{ {{\zeta ^ \star }\left( {p,m - 1} \right) + {\zeta ^ \star }\left( {m,p - 1} \right) - \zeta \left( {m + p - 1} \right) - \zeta \left( m \right)\zeta \left( p \right)} \right\}\nonumber\\
& + \left\{ {{\zeta ^ \star }\left( {m + p - 1,r} \right) - {\zeta ^ \star }\left( {m + p,r - 1} \right)} \right\}\nonumber\\
& - \left\{ {{\zeta ^ \star }\left( {p,m - 1,r} \right) - {\zeta ^ \star }\left( {p,m,r - 1} \right)} \right\}\nonumber\\
& - \left\{ {{\zeta ^ \star }\left( {m,p - 1,r} \right) - {\zeta ^ \star }\left( {m,p,r - 1} \right)} \right\}.
\end{align}
By the definition of multiple zeta function and multiple zeta-star function, the following relation is easily derived
\begin{align}\label{3.7}
\zeta^\star\left( {q,p,r} \right) = \zeta \left( {q,p,r} \right) + \zeta\left( {q,p + r} \right) + \zeta\left( {q + p,r} \right) + \zeta \left( {q +p+ r} \right).
\end{align}
Hence, by a direct calculation, we can rewrite (\ref{3.7}) as
\begin{align}\label{3.8}
F\left( {\zeta \left( m \right),\zeta \left( p \right),\zeta \left( r \right)} \right) =&\zeta \left( {p,m,r - 1} \right) + \zeta \left( {p,r,m - 1} \right) + \zeta \left( {m,p,r - 1} \right) + \zeta \left( {m,r,p - 1} \right)\nonumber\\
& + \zeta \left( {r,p,m - 1} \right) + \zeta \left( {r,m,p - 1} \right) + \zeta \left( {m + p,r - 1} \right) + \zeta \left( {p + r,m - 1} \right)\nonumber\\
&+ \zeta \left( {m + r,p - 1} \right) + \zeta \left( {p,m + r - 1} \right) + \zeta \left( {m,p + r - 1} \right) \nonumber\\
& + \zeta \left( {r,p + m - 1} \right)+ \zeta \left( {m + p + r - 1} \right) - \zeta \left( p \right)\zeta \left( m \right)\zeta \left( r \right).
\end{align}
This formula was also given by Hoffman in \cite{MEH2016}.

In (\ref{1.11}), letting $n\rightarrow \infty$, then multiplying it by $\frac{{{{\ln }^{r - 1}}x}}{{x\left( {1 - x} \right)}}\quad (r\in \N\setminus\{1\})$
and integrating over the interval (0,1), we conclude that
\begin{align}\label{3.9}
\sum\limits_{n = 1}^\infty  {\frac{{\left( {A - {A_n}} \right)\left( {B - {B_n}} \right)}}{{{n^r}}}}  =& \sum\limits_{n = 1}^\infty  {\left( {\sum\limits_{k = 1}^n {{a_k}{\zeta _{k - 1}}\left( r \right)} } \right){b_n}}  + \sum\limits_{n = 1}^\infty  {\left( {\sum\limits_{k = 1}^n {{b_k}{\zeta _{k - 1}}\left( r \right)} } \right){a_n}}\nonumber\\
&  - \sum\limits_{n = 1}^\infty  {{a_n}{b_n}{\zeta _{n - 1}}\left( r \right)}.
\end{align}
In fact, it is easily shown that the formula (\ref{3.9}) is also true when $r=1$.\\
Taking $\left( {{a_n},{b_n}} \right) = \left( {\frac{1}{{{n^m}}},\frac{1}{{{n^p}}}} \right)$ in (\ref{3.9}) and using the relation (\ref{3.7}), we get
\begin{align}\label{3.10}
\sum\limits_{n = 1}^\infty  {\frac{{\left( {\zeta \left( m \right) - {\zeta _n}\left( m \right)} \right)\left( {\zeta \left( p \right) - {\zeta _n}\left( p \right)} \right)}}{{{n^r}}}} =\zeta \left( {m,p,r} \right) + \zeta \left( {p,m,r} \right) + \zeta \left( {m + p,r} \right).
\end{align}
From definition of the multiple harmonic star sum, it is easily shown that
\[\zeta_n^\star\left(p,m\right)+\zeta_n^\star\left(m,p\right)=\zeta_n\left(p\right)\zeta_n\left(m\right)+\zeta_n\left(p+m\right).\]
Hence, we give the general formula
\[\sum\limits_{n = 1}^\infty  {\frac{{{\zeta _n}\left( m \right){\zeta _n}\left( p \right)}}{{{n^r}}}}  = {\zeta ^ \star }\left( {r,m,p} \right) + {\zeta ^ \star }\left( {r,p,m} \right) - {\zeta ^ \star }\left( {r,m + p} \right)\quad (r>1,m,p\geq1).\]
Moreover, by a simple calculation, we derive the following identity
\begin{align}\label{3.11}
\sum\limits_{n = 1}^\infty  {\frac{{{\zeta _n}\left( p \right){\zeta _n}\left( r \right)}}{{{n^m}}}}  = {\zeta ^ \star }\left( {p + m,r} \right) + \zeta \left( p \right){\zeta ^ \star }\left( {m,r} \right)-{\zeta ^ \star }\left( {p,m,r} \right),\ p,m\in \N\setminus\{1\}, r\in \N.
\end{align}
Substituting (\ref{3.11}) into (\ref{3.10}), we obtain the symmetry identity
\begin{align}\label{3.12}
&\sum\limits_{n = 1}^\infty  {\left\{ {\frac{{{\zeta _n}\left( m \right){\zeta _n}\left( p \right)}}{{{n^r}}} + \frac{{{\zeta _n}\left( p \right){\zeta _n}\left( r \right)}}{{{n^m}}} + \frac{{{\zeta _n}\left( m \right){\zeta _n}\left( r \right)}}{{{n^p}}}} \right\}}\nonumber \\
&={\zeta ^ \star }\left( {p + m,r} \right) + {\zeta ^ \star }\left( {p + r,m} \right) + {\zeta ^ \star }\left( {m + r,p} \right)\nonumber\\
&\quad + \zeta \left( p \right)\zeta \left( m \right)\zeta \left( r \right) - \zeta \left( {p + m + r} \right),\ \ \ \ \ \ \ \ \ \ \ \ \  p,m,r\in \N/\{1\}.
\end{align}
From formula (\ref{3.11}) and references \cite{BG1996,FS1998,X2016,X2017,Xu2017}, we know that the quadratic sums
 \[S_{mp,r}:=\sum\limits_{n = 1}^\infty  {\frac{{{\zeta _n}\left( m \right){\zeta _n}\left( p \right)}}{{{n^r}}}} \]
can be evaluated in terms of zeta values and double zeta (star) values in the following cases: $m=p=1$, $m=p=r$ and $m+p+r$ is even or less than or equal to $m+p+r\leq 10$. For example, from \cite{X2016,X2017} we have
\begin{align*}
&\sum\limits_{n = 1}^\infty  {\frac{{{H_n}{\zeta _n}\left( 3 \right)}}{{{n^3}}}}  = \frac{{83}}{8}\zeta \left( 7 \right) + \frac{1}{4}\zeta \left( 3 \right)\zeta \left( 4 \right) - \frac{{11}}{2}\zeta \left( 2 \right)\zeta \left( 5 \right),\\
& \sum\limits_{n = 1}^\infty  {\frac{{{H_n}{\zeta _n}\left( 2 \right)}}{{{n^5}}}}  =  - \frac{{343}}{{48}}\zeta \left( 8 \right) + 12\zeta \left( 3 \right)\zeta \left( 5 \right) - \frac{5}{2}\zeta \left( 2 \right){\zeta ^2}\left( 3 \right) - \frac{3}{4}\zeta^\star(6,2),\\
&\sum\limits_{n = 1}^\infty  {\frac{{\zeta _n^2\left( 2 \right)}}{{{n^4}}}}  = 11\zeta^\star(6,2) + \frac{{457}}{{18}}\zeta \left( 8 \right) + 6\zeta \left( 2 \right){\zeta ^2}\left( 3 \right) - 40\zeta \left( 3 \right)\zeta \left( 5 \right),\\
&\sum\limits_{n = 1}^\infty  {\frac{{\zeta _n^2\left( 2 \right){\zeta _n}\left( 3 \right)}}
{{{n^2}}}}  =  - \frac{{617}}
{{72}}\zeta (9) + {\zeta ^3}(3) + \frac{{91}}
{8}\zeta (2)\zeta (7) - \frac{{17}}
{4}\zeta (4)\zeta (5) - \frac{{329}}
{{84}}\zeta (3)\zeta (6),\\
&\sum\limits_{n = 1}^\infty  {\frac{{\zeta _n^2\left( 2 \right)}}
{{{n^6}}}}  = \frac{{2697}}
{{40}}\zeta (10) - 41{\zeta ^2}(5) - 63\zeta (3)\zeta (7) + 16\zeta (2)\zeta (3)\zeta (5) + 4{\zeta ^2}(3)\zeta (4),\\
&\quad\quad\quad\quad\quad\quad + \frac{{23}}
{2}{\zeta ^ \star }(8,2){\text{  +  }}2\zeta (2){\zeta ^ \star }(6,2).
\end{align*}
Similarly, from (\ref{1.11}), and by a similar as method in the proofs of (\ref{3.3})-(\ref{3.5}), we deduce the following results
\begin{align}\label{3.13}
&\sum\limits_{n = 1}^\infty  {\left( {\zeta \left( m \right) - {\zeta _n}\left( m \right)} \right)\left( {\zeta \left( p \right) - {\zeta _n}\left( p \right)} \right){{\left( { - 1} \right)}^{n - 1}}}\nonumber \\
& = \frac{1}{2}\left\{ {\zeta \left( m \right)\zeta \left( p \right) + \zeta \left( {m,\bar p} \right) + \zeta \left( {p,\bar m} \right) + \zeta \left( \overline{m + p} \right)} \right\},
\end{align}
\begin{align}\label{3.14}
&\sum\limits_{n = 1}^\infty  {\left( {{\zeta ^ \star }\left( {p + 1,p} \right) - \zeta _n^ \star \left( {p + 1,p} \right)} \right)\left( {{\zeta ^ \star }\left( {q + 1,q} \right) - \zeta _n^ \star \left( {q + 1,q} \right)} \right)} \nonumber \\
& = \frac{1}{2}\sum\limits_{n = 1}^\infty  {\left\{ {\frac{{\zeta _n^2\left( p \right){\zeta _n}\left( q \right) + {\zeta _n}\left( {2p} \right){\zeta _n}\left( q \right)}}{{{n^{q + 1}}}} + \frac{{\zeta _n^2\left( q \right){\zeta _n}\left( p \right) + {\zeta _n}\left( {2q} \right){\zeta _n}\left( p \right)}}{{{n^{p + 1}}}}} \right\}}\nonumber \\
&\quad - \sum\limits_{n = 1}^\infty  {\frac{{{\zeta _n}\left( p \right){\zeta _n}\left( q \right)}}{{{n^{p + q + 1}}}}}  - \left( {\sum\limits_{n = 1}^\infty  {\frac{{{\zeta _n}\left( p \right)}}{{{n^{p + 1}}}}} } \right)\left( {\sum\limits_{n = 1}^\infty  {\frac{{{\zeta _n}\left( q \right)}}{{{n^{q + 1}}}}} } \right),
\end{align}
\begin{align}\label{3.15}
&\sum\limits_{n = 1}^\infty  {\left( {{\zeta ^ \star }\left( {p + 1,p} \right) - \zeta _n^ \star \left( {p + 1,p} \right)} \right)\left( {{\zeta ^ \star }\left( {p + 2,p} \right) - \zeta _n^ \star \left( {p + 2,p} \right)} \right)} \nonumber \\
& = \frac{1}{2}\sum\limits_{n = 1}^\infty  {\left\{ {\frac{{\zeta _n^3\left( p \right) + {\zeta _n}\left( p \right){\zeta _n}\left( {2p} \right)}}{{{n^{p + 2}}}} - \frac{{\zeta _n^2\left( p \right)}}{{{n^{2p + 2}}}}} \right\}}  + \frac{1}{2}{\left( {\sum\limits_{n = 1}^\infty  {\frac{{{\zeta _n}\left( p \right)}}{{{n^{p + 1}}}}} } \right)^2}\nonumber\\
& \quad- \left( {\sum\limits_{n = 1}^\infty  {\frac{{{\zeta _n}\left( p \right)}}{{{n^{p + 1}}}}} } \right)\left( {\sum\limits_{n = 1}^\infty  {\frac{{{\zeta _n}\left( p \right)}}{{{n^{p + 2}}}}} } \right).
\end{align}
Next, we give a Theorem.
\begin{thm}\label{thm3.1}
Let $m,p\ge 2$ be integers, then we have
\begin{align}\label{3.16}
\begin{array}{*{20}{c}}
   {\sum\limits_{n = 1}^\infty  {\frac{{\zeta \left( m \right){\zeta _n}\left( p \right) - \zeta \left( p \right){\zeta _n}\left( m \right)}}{n}}  = \zeta \left( p \right)\zeta \left( {m,1} \right) - \zeta \left( m \right)} \hfill  \\
\end{array}\zeta \left( {p,1} \right).
\end{align}
\end{thm}
Next, we will use two different methods to prove the theorem above.\\
\pf{\bf 1.} First, we give a complicated proof.  To prove the identity (\ref{3.16}), we construct the following power function
\[y := \sum\limits_{n = 1}^\infty  {\left\{ {{H_n}{\zeta _n}\left( m \right) - {\zeta _n}\left( {m + 1} \right)} \right\}{x^{n - 1}}},\ x\in (-1,1). \]
Then, by using the definition of harmonic numbers, it is easily seen that
\begin{align}\label{3.17}
y &= \sum\limits_{n = 1}^\infty  {\left\{ {\left( {{H_n} + \frac{1}{{n + 1}}} \right)\left( {{\zeta _n}\left( m \right) + \frac{1}{{{{\left( {n + 1} \right)}^m}}}} \right) - \left( {{\zeta _n}\left( {m + 1} \right) + \frac{1}{{{{\left( {n + 1} \right)}^{m + 1}}}}} \right)} \right\}{x^n}}
\nonumber \\ &=\sum\limits_{n = 1}^\infty  {\left\{ {{H_n}{\zeta _n}\left( m \right) - {\zeta _n}\left( {m + 1} \right) + \frac{{{H_n}}}{{{{\left( {n + 1} \right)}^m}}} + \frac{{{\zeta _n}\left( m \right)}}{{n + 1}}} \right\}{x^n}} .
\end{align}
Therefore, from (\ref{3.17}), we obtain the following relation
\begin{align}\label{3.18}
\sum\limits_{n = 1}^\infty  {\left\{ {{H_n}{\zeta _n}\left( m \right) - {\zeta _n}\left( {m + 1} \right)} \right\}{x^{n - 1}}}  = \sum\limits_{n = 1}^\infty  {\left\{ {\frac{{{H_n}}}{{{{\left( {n + 1} \right)}^m}}} + \frac{{{\zeta _n}\left( m \right)}}{{n + 1}}} \right\}\frac{{{x^n}}}{{1 - x}}}.
\end{align}
Then multiplying (\ref{3.18}) by ${\ln ^{p - 1}}x$ and integrating over (0,1), we conclude that
\begin{align}\label{3.19}
\sum\limits_{n = 1}^\infty  {\frac{{{H_n}{\zeta _n}\left( m \right) - {\zeta _n}\left( {m + 1} \right)}}{{{n^p}}}}  = \sum\limits_{n = 1}^\infty  {\left\{ {\frac{{{H_n}}}{{{{\left( {n + 1} \right)}^m}}} + \frac{{{\zeta _n}\left( m \right)}}{{n + 1}}} \right\}\left\{ {\zeta \left( p \right) - {\zeta _n}\left( p \right)} \right\}}.
\end{align}
Hence, by a direct calculation, the result is
\begin{align}\label{3.20}
\sum\limits_{n = 1}^\infty  {\left\{ {\frac{{{H_n}{\zeta _n}\left( m \right)}}{{{n^p}}} + \frac{{{H_n}{\zeta _n}\left( p \right)}}{{{n^m}}}} \right\}} &=\zeta \left( p \right)\sum\limits_{n = 1}^\infty  {\frac{{{H_n}}}{{{n^m}}}}  + \sum\limits_{n = 1}^\infty  {\frac{{{H_n}}}{{{n^{p + m}}}}}  + \sum\limits_{n = 1}^\infty  {\frac{{{\zeta _n}\left( m \right)}}{{{n^{p + 1}}}}}
\nonumber \\ &\quad- \sum\limits_{n = 1}^\infty  {\frac{{{\zeta _n}\left( {m + 1} \right)}}{{{n^p}}}}  + \sum\limits_{n = 1}^\infty  {\frac{{{\zeta _n}\left( m \right)}}{n}} \left\{ {\zeta \left( p \right) - {\zeta _n}\left( p \right)} \right\}.
\end{align}
Changing $m,p$ to $p,m$ in above equation, we can get
\begin{align}\label{3.21}
\sum\limits_{n = 1}^\infty  {\left\{ {\frac{{{H_n}{\zeta _n}\left( p \right)}}{{{n^m}}} + \frac{{{H_n}{\zeta _n}\left( m \right)}}{{{n^p}}}} \right\}} &=\zeta \left( m \right)\sum\limits_{n = 1}^\infty  {\frac{{{H_n}}}{{{n^p}}}}  + \sum\limits_{n = 1}^\infty  {\frac{{{H_n}}}{{{n^{m + p}}}}}  + \sum\limits_{n = 1}^\infty  {\frac{{{\zeta _n}\left( p \right)}}{{{n^{m + 1}}}}}
\nonumber \\ &\quad - \sum\limits_{n = 1}^\infty  {\frac{{{\zeta _n}\left( {p + 1} \right)}}{{{n^m}}}}  + \sum\limits_{n = 1}^\infty  {\frac{{{\zeta _n}\left( p \right)}}{n}} \left\{ {\zeta \left( m \right) - {\zeta _n}\left( m \right)} \right\}.
\end{align}
Thus, combining (\ref{3.20}) with (\ref{3.21}), and using the simple relation
\[\sum\limits_{n = 1}^\infty  {\frac{{{H_n}}}{{{n^m}}}}  = \zeta \left( {m,1} \right) + \zeta \left( {m + 1} \right),\]
we may deduce the desired result. The proof of Theorem \ref{thm3.1} is finished.\hfill$\square$\\
\pf {\bf 2.} In fact, one can prove a more general result without using the above complicated way. Writing $\zeta_{>n}$ for $\zeta(a)-\zeta_n(a)$, one has
\begin{align*}
   &{\sum\limits_{n = 1}^\infty  {\frac{{\zeta \left( m \right){\zeta _n}\left( p \right) - \zeta \left( p \right){\zeta _n}\left( m \right)}}{{{n^k}}}} } \hfill  \\
  &= \sum\limits_{n = 1}^\infty  {\frac{{\left( {{\zeta _n}\left( m \right) + {\zeta _{ > n}}\left( m \right)} \right){\zeta _n}\left( p \right) - \left( {{\zeta _n}\left( p \right) + {\zeta _{ > n}}\left( p \right)} \right){\zeta _n}\left( m \right)}}{{{n^k}}}}  \\
  &= \sum\limits_{n = 1}^\infty  {\frac{{{\zeta _{ > n}}\left( m \right){\zeta _n}\left( p \right) - {\zeta _{ > n}}\left( p \right){\zeta _n}\left( m \right)}}{{{n^k}}}}  \\
  &= \zeta \left( {m,k,p} \right) + \zeta \left( {m,k + p} \right) - \zeta \left( {p,k,m} \right) - \zeta \left( {p,k + m} \right)
\end{align*}
for any positive integer $k$. One the other hand, by the well-known ``harmonic product" (also called the ``stuffle product") of multiple zeta values \cite{MEH1997}, we have
\[\zeta \left( p \right)\zeta \left( {m,k} \right) = \zeta \left( {p,m,k} \right) + \zeta \left( {p + m,k} \right) + \zeta \left( {m,p,k} \right) + \zeta \left( {m,p + k} \right) + \zeta \left( {m,k,p} \right).\]
Similarly expanding out $\zeta \left( m \right)\zeta \left( {p,k} \right)$ and cancelling gives
\begin{align}\label{3.22}
\begin{array}{*{20}{c}}
   {\sum\limits_{n = 1}^\infty  {\frac{{\zeta \left( m \right){\zeta _n}\left( p \right) - \zeta \left( p \right){\zeta _n}\left( m \right)}}{{{n^k}}}} } \hfill  \\
\end{array} = \zeta \left( p \right)\zeta \left( {m,k} \right) - \zeta \left( m \right)\zeta \left( {p,k} \right).\end{align}
It is clear that (\ref{3.16}) is immediate corollary of (\ref{3.22}).\hfill$\square$

Next, we continue to establish some relations between sums of zeta tails and multiple zeta (star) values.
\begin{align}\label{3.23}
&\sum\limits_{n = 1}^\infty  {\frac{{\zeta \left( m \right)\zeta \left( p \right) - {\zeta _n}\left( m \right){\zeta _n}\left( p \right)}}{n}}\nonumber  \\&= \zeta \left( p \right)\sum\limits_{n = 1}^\infty  {\frac{{\zeta \left( m \right) - {\zeta _n}\left( m \right)}}{n}}  + \sum\limits_{n = 1}^\infty  {\frac{{{\zeta _n}\left( m \right)\left( {\zeta \left( p \right) - {\zeta _n}\left( p \right)} \right)}}{n}},\ m,p\in \N\setminus\{1\}.
\end{align}
In \cite{OF2016}, O. Furdui and C. V\u{a}lean gave the formula
\begin{align}\label{3.24}
\sum\limits_{n = 1}^\infty  {\frac{{\zeta \left( m \right) - {\zeta _n}\left( m \right)}}{n}}  = {\zeta}\left( {m,1} \right) .
\end{align}
Hence, combining (\ref{3.20}), (\ref{3.23}) and (\ref{3.24}), we have the following corollary.
\begin{cor}\label{cor3.2} For positive integers $m>1$ and $p>1$, then we have
\begin{align}\label{3.25}
&\sum\limits_{n = 1}^\infty  {\frac{{\zeta \left( m \right)\zeta \left( p \right) - {\zeta _n}\left( m \right){\zeta _n}\left( p \right)}}{n}}\nonumber \\
& = \sum\limits_{n = 1}^\infty  {\left\{ {\frac{{{H_n}{\zeta _n}\left( m \right)}}{{{n^p}}} + \frac{{{H_n}{\zeta _n}\left( p \right)}}{{{n^m}}}} \right\}}  + {\zeta ^ \star }\left( {p,m + 1} \right)\nonumber\\
&\quad - \zeta \left( p \right)\zeta \left( {m + 1} \right) - {\zeta ^ \star }\left( {p + m,1} \right) - {\zeta ^ \star }\left( {p + 1,m} \right).
\end{align}
\end{cor}
Similarly, noting that
\begin{align*}
&\sum\limits_{n = 1}^\infty  {\frac{{\zeta_n \left( m \right)\zeta \left( p \right)\zeta \left( q \right) - \zeta \left( m \right){\zeta _n}\left( p \right){\zeta _n}\left( q \right)}}{n}} \\
& = \zeta \left( m \right)\sum\limits_{n = 1}^\infty  {\frac{{\zeta \left( p \right)\zeta \left( q \right) - {\zeta _n}\left( p \right){\zeta _n}\left( q \right)}}{n}} - \zeta \left( p \right)\zeta \left( q \right)\sum\limits_{n = 1}^\infty  {\frac{{\zeta \left( m \right) - {\zeta _n}\left( m \right)}}{n}},\ \ m,p,q\in \N\setminus\{1\} ,
\end{align*}
 we obtain the following corollary.
\begin{cor}\label{cor3.3}For positive integers $m,p,q>1$, then we have
\begin{align}\label{3.26}
&\sum\limits_{n = 1}^\infty  {\frac{{{\zeta _n}\left( m \right)\zeta \left( p \right)\zeta \left( q \right) - \zeta \left( m \right){\zeta _n}\left( p \right){\zeta _n}\left( q \right)}}{n}} \nonumber\\
& = \zeta \left( m \right)\left\{ \begin{array}{l}
 \sum\limits_{n = 1}^\infty  {\left( {\frac{{{H_n}{\zeta _n}\left( q \right)}}{{{n^p}}} + \frac{{{H_n}{\zeta _n}\left( p \right)}}{{{n^q}}}} \right)}  + {\zeta ^ \star }\left( {p,q + 1} \right)\nonumber \\
  - {\zeta ^ \star }\left( {p + q,1} \right) - {\zeta ^ \star }\left( {p + 1,q} \right) - \zeta \left( p \right)\zeta \left( {q + 1} \right) \\
 \end{array} \right\}\\
&\quad - \zeta \left( p \right)\zeta \left( q \right){{\zeta }\left( {m,1} \right) } .
\end{align}
\end{cor}
Furthermore, by considering the following power series function
\[y = \sum\limits_{n = 1}^\infty  {\left\{ {{H_n}{\zeta _n}\left( p \right){\zeta _n}\left( q \right) - {\zeta _n}\left( {p + q + 1} \right)} \right\}{x^{n - 1}}} ,\;x \in ( - 1,1),\]
and proceeding in a similar fashion to evaluation of the Theorem \ref{thm3.1}, we can give the following
Theorem.
\begin{thm}\label{thm3.4} For positive integers $p,m>1$ and $q>0$, we have
\begin{align*}
&\sum\limits_{n = 1}^\infty  {\frac{{\zeta \left( m \right){\zeta _n}\left( p \right){\zeta _n}\left( q \right) - \zeta \left( p \right){\zeta _n}\left( m \right){\zeta _n}\left( q \right)}}{n}} \\
& = \zeta \left( p \right)\left\{ {\sum\limits_{n = 1}^\infty  {\left( {\frac{{{H_n}{\zeta _n}\left( m \right)}}{{{n^q}}} + \frac{{{H_n}{\zeta _n}\left( q \right)}}{{{n^m}}}} \right)}  - {\zeta ^ \star }\left( {q + 1,m} \right) - {\zeta ^ \star }\left( {m + q,1} \right) - {\zeta ^ \star }\left( {m + 1,q} \right)} \right\}\\
&\quad - \zeta \left( m \right)\left\{ {\sum\limits_{n = 1}^\infty  {\left( {\frac{{{H_n}{\zeta _n}\left( q \right)}}{{{n^p}}} + \frac{{{H_n}{\zeta _n}\left( p \right)}}{{{n^q}}}} \right)}  - {\zeta ^ \star }\left( {q + 1,p} \right) - {\zeta ^ \star }\left( {p + q,1} \right) - {\zeta ^ \star }\left( {p + 1,q} \right)} \right\}\\
&\quad + \zeta \left( {m + q + 1} \right)\zeta \left( p \right) - \zeta \left( {p + q + 1} \right)\zeta \left( m \right).
\end{align*}
\end{thm}
On the other hand, we note that
\begin{align}\label{3.27}
\sum\limits_{n = 1}^\infty  {\frac{{{\zeta _n}\left( m \right)\left( {\zeta \left( p \right) - {\zeta _n}\left( p \right)} \right)}}{n}}  = \zeta \left( {p,1,m} \right) + \zeta \left( {p,m + 1} \right).
\end{align}
The relations (\ref{3.20}) and (\ref{3.27}) yield the following result:
\begin{align}\label{3.28}
 \sum\limits_{n = 1}^\infty  {\left\{ {\frac{{{H_n}{\zeta _n}\left( m \right)}}{{{n^p}}} + \frac{{{H_n}{\zeta _n}\left( p \right)}}{{{n^m}}}} \right\}}  =& \zeta \left( p \right)\zeta \left( {m,1} \right) + \zeta \left( p \right)\zeta \left( {m + 1} \right) + \zeta \left( {p + m,1} \right)\nonumber \\
  &+ \zeta \left( {p + 1,m} \right) + \zeta \left( {p + m + 1} \right) + \zeta \left( {p,1,m} \right).
\end{align}
Therefore, by using the above relation (\ref{3.28}), Corollary \ref{cor3.2}, Corollary \ref{cor3.3} and Theorem \ref{thm3.4} can be written as follows.
\begin{cor}\label{cor3.5} For positive integers $m>1$ and $p>1$, then the following relations hold:
\begin{align}\label{3.29}
&\sum\limits_{n = 1}^\infty  {\frac{{\zeta \left( m \right)\zeta \left( p \right) - {\zeta _n}\left( m \right){\zeta _n}\left( p \right)}}{n}}\nonumber \\
& =\zeta \left( p \right)\zeta \left( {m,1} \right) + \zeta \left( {p,1,m} \right) + \zeta \left( {p,m + 1} \right).
\end{align}
\begin{align}\label{3.30}
&\sum\limits_{n = 1}^\infty  {\frac{{{\zeta _n}\left( m \right)\zeta \left( p \right)\zeta \left( q \right) - \zeta \left( m \right){\zeta _n}\left( p \right){\zeta _n}\left( q \right)}}{n}} \nonumber\\
&=\zeta \left( m \right)\left[ {\zeta \left( p \right)\zeta \left( {q,1} \right) + \zeta \left( {p,1,q} \right) + \zeta \left( {p,q + 1} \right)} \right] - \zeta \left( p \right)\zeta \left( q \right)\zeta \left( {m,1} \right),
\end{align}
\begin{align}\label{3.31}
&\sum\limits_{n = 1}^\infty  {\frac{{\zeta \left( m \right){\zeta _n}\left( p \right){\zeta _n}\left( q \right) - \zeta \left( p \right){\zeta _n}\left( m \right){\zeta _n}\left( q \right)}}{n}} \nonumber\\
&=\zeta \left( p \right)\left[ {\zeta \left( {m,1,q} \right) + \zeta \left( {m,q + 1} \right)} \right] - \zeta \left( m \right)\left[ {\zeta \left( {p,1,q} \right) + \zeta \left( {p,q + 1} \right)} \right].
\end{align}
\end{cor}
In fact, applying the same arguments as in the proof of (\ref{3.22}), the above formulas (\ref{3.29})-(\ref{3.31}) can be also proved by using the method of stuffle product.

We now close this section with a final theorem.
\begin{thm}
For positive integers $m,p,k>1$, then the following relations hold:
\begin{align}\label{3.32}
&\sum\limits_{n = 1}^\infty  {\frac{{\zeta \left( m \right)\zeta \left( p \right) - {\zeta _n}\left( m \right){\zeta _n}\left( p \right)}}{n^k}}\nonumber \\
& =\zeta \left( p \right)\zeta \left( {m,k} \right) + \zeta \left( {p,k,m} \right) + \zeta \left( {p,m + k} \right).
\end{align}
\begin{align}\label{3.33}
&\sum\limits_{n = 1}^\infty  {\frac{{{\zeta _n}\left( m \right)\zeta \left( p \right)\zeta \left( q \right) - \zeta \left( m \right){\zeta _n}\left( p \right){\zeta _n}\left( q \right)}}{n^k}} \nonumber\\
&=\zeta \left( m \right)\left[ {\zeta \left( p \right)\zeta \left( {q,k} \right) + \zeta \left( {p,k,q} \right) + \zeta \left( {p,q + k} \right)} \right] - \zeta \left( p \right)\zeta \left( q \right)\zeta \left( {m,k} \right),
\end{align}
\begin{align}\label{3.34}
&\sum\limits_{n = 1}^\infty  {\frac{{\zeta \left( m \right){\zeta _n}\left( p \right){\zeta _n}\left( q \right) - \zeta \left( p \right){\zeta _n}\left( m \right){\zeta _n}\left( q \right)}}{n^k}} \nonumber\\
&=\zeta \left( p \right)\left[ {\zeta \left( {m,k,q} \right) + \zeta \left( {m,q + k} \right)} \right] - \zeta \left( m \right)\left[ {\zeta \left( {p,k,q} \right) + \zeta \left( {p,q + k} \right)} \right].
\end{align}
\end{thm}
\pf On the one hand, for any integers $m,p,k>1$, we note that the series on the left hand side of formulas (\ref{3.32})-(\ref{3.34}) can be written as
\begin{align}\label{3.35}
\sum\limits_{n = 1}^\infty  {\frac{{\zeta \left( m \right)\zeta \left( p \right) - {\zeta _n}\left( m \right){\zeta _n}\left( p \right)}}{{{n^k}}}}  = \zeta \left( m \right)\zeta \left( p \right)\zeta \left( k \right) - \sum\limits_{n = 1}^\infty  {\frac{{{\zeta _n}\left( m \right){\zeta _n}\left( p \right)}}{{{n^k}}}} ,
\end{align}
\begin{align}\label{3.36}
\sum\limits_{n = 1}^\infty  {\frac{{{\zeta _n}\left( m \right)\zeta \left( p \right)\zeta \left( q \right) - \zeta \left( m \right){\zeta _n}\left( p \right){\zeta _n}\left( q \right)}}{{{n^k}}}}  = \zeta \left( p \right)\zeta \left( q \right){\zeta ^ \star }\left( {k,m} \right) - \zeta \left( m \right)\sum\limits_{n = 1}^\infty  {\frac{{{\zeta _n}\left( p \right){\zeta _n}\left( q \right)}}{{{n^k}}}} ,
\end{align}
\begin{align}\label{3.37}
\sum\limits_{n = 1}^\infty  {\frac{{\zeta \left( m \right){\zeta _n}\left( p \right){\zeta _n}\left( q \right) - \zeta \left( p \right){\zeta _n}\left( m \right){\zeta _n}\left( q \right)}}{{{n^k}}}}  = \zeta \left( m \right)\sum\limits_{n = 1}^\infty  {\frac{{{\zeta _n}\left( p \right){\zeta _n}\left( q \right)}}{{{n^k}}}}  - \zeta \left( p \right)\sum\limits_{n = 1}^\infty  {\frac{{{\zeta _n}\left( m \right){\zeta _n}\left( q \right)}}{{{n^k}}}}.
\end{align}
On the other hand, from (\ref{3.7}) and (\ref{3.11}), we readily find that
\begin{align}\label{3.38}
\sum\limits_{n = 1}^\infty  {\frac{{{\zeta _n}\left( m \right){\zeta _n}\left( p \right)}}{{{n^k}}}}  = \zeta \left( p \right)\zeta \left( m \right)\zeta \left( k \right) - \zeta \left( p \right)\zeta \left( {m,k} \right) - \zeta \left( {p,k,m} \right) - \zeta \left( {p,k + m} \right).
\end{align}
Then, substituting (\ref{3.38}) into (\ref{3.35})-(\ref{3.37}), by a simple calculation, we may easily deduce the desired results. \hfill$\square$\\
It is clear that the identities (\ref{3.32})-(\ref{3.34}) can be proved by using the method of stuffle product.
\section{Some Examples}
In section 3, we prove many interesting results involving zeta tails. In this section, we apply these results to give various specific cases.
Some illustrative examples follow.
\begin{align*}
&\sum\limits_{n = 1}^\infty  {\frac{{{{\left( {\zeta \left( 2 \right) - {\zeta _n}\left( 2 \right)} \right)}^2}}}{n}}  = 5\zeta \left( 2 \right)\zeta \left( 3 \right) - 9\zeta \left( 5 \right),\\
&\sum\limits_{n = 1}^\infty  {{{\left( {\bar \zeta \left( 2 \right) - {{\bar \zeta }_n}\left( 2 \right)} \right)}^2}}  = 3\zeta \left( 2 \right)\ln 2 - \frac{9}{4}\zeta \left( 3 \right) - \frac{5}{8}\zeta \left( 4 \right),\\
&\sum\limits_{n = 1}^\infty  {\frac{{\left( {\zeta \left( 2 \right) - {\zeta _n}\left( 2 \right)} \right)\left( {\zeta \left( 3 \right) - {\zeta _n}\left( 3 \right)} \right)}}{n}}  = {\zeta ^2}\left( 3 \right) - \frac{{61}}{{48}}\zeta \left( 6 \right),\\
&\sum\limits_{n = 1}^\infty  {{{\left( {\zeta \left( 2 \right) - {\zeta _n}\left( 2 \right)} \right)}^3}}  = 9\zeta \left( 2 \right)\zeta \left( 3 \right) - \frac{{35}}{8}\zeta \left( 6 \right) - \frac{{25}}{2}\zeta \left( 5 \right),\\
&\sum\limits_{n = 1}^\infty  {{{\left( {{\zeta ^ \star }\left( {2,1} \right) - \zeta _n^ \star \left( {2,1} \right)} \right)}^2}}  = \frac{{15}}{2}\zeta \left( 5 \right) + 3\zeta \left( 2 \right)\zeta \left( 3 \right) - 4{\zeta ^2}\left( 3 \right),\\
&\sum\limits_{n = 1}^\infty  {\frac{{{\zeta ^2}\left( 3 \right) - \zeta _n^2\left( 3 \right)}}{n}}  =  - \frac{{73}}{4}\zeta \left( 7 \right) + \frac{3}{2}\zeta \left( 3 \right)\zeta \left( 4 \right) + 10\zeta \left( 2 \right)\zeta \left( 5 \right),\\
&\sum\limits_{n = 1}^\infty  {\left( {\zeta \left( 2 \right) - {\zeta _n}\left( 2 \right)} \right)\left( {\bar \zeta \left( 2 \right) - {{\bar \zeta }_n}\left( 2 \right)} \right)}  = \frac{3}{2}\zeta \left( 2 \right)\ln 2 - \frac{5}{4}\zeta \left( 4 \right) - \frac{3}{8}\zeta \left( 3 \right),\\
&\sum\limits_{n = 1}^\infty  {\left( {\zeta \left( 2 \right) - {\zeta _n}\left( 2 \right)} \right)\left( {\zeta \left( 3 \right) - {\zeta _n}\left( 3 \right)} \right){{\left( { - 1} \right)}^{n - 1}}}  = \frac{9}{{16}}\zeta \left( 2 \right)\zeta \left( 3 \right) - \frac{{31}}{{32}}\zeta \left( 5 \right),\\
&\sum\limits_{n = 1}^\infty  {{{\left( {\zeta \left( 2 \right) - {\zeta _n}\left( 2 \right)} \right)}^2}} \left( {\zeta \left( 3 \right) - {\zeta _n}\left( 3 \right)} \right) = \frac{7}{6}\zeta \left( 6 \right) + \frac{3}{2}{\zeta ^2}\left( 3 \right) - \frac{5}{2}\zeta \left( 3 \right)\zeta \left( 4 \right),\\
&\sum\limits_{n = 1}^\infty  {\left( {{\zeta ^ \star }\left( {2,1} \right) - \zeta _n^ \star \left( {2,1} \right)} \right)\left( {{\zeta ^ \star }\left( {3,1} \right) - \zeta _n^ \star \left( {3,1} \right)} \right)}  =  - \frac{1}{3}\zeta \left( 6 \right) + 3{\zeta ^2}\left( 3 \right) - \frac{5}{2}\zeta \left( 3 \right)\zeta \left( 4 \right),\\
&\sum\limits_{n = 1}^\infty  {{{\left( {\zeta \left( 2 \right) - {\zeta _n}\left( 2 \right)} \right)}^2}{{\left( { - 1} \right)}^{n - 1}}}  = 4{\rm{L}}{{\rm{i}}_4}\left( {\frac{1}{2}} \right) - \frac{{29}}{8}\zeta \left( 4 \right) - \zeta \left( 2 \right){\ln ^2}2 + \frac{1}{6}{\ln ^4}2 + \frac{7}{2}\zeta \left( 3 \right)\ln 2,\\
&\sum\limits_{n = 1}^\infty  {\frac{{{\zeta _n}\left( 2 \right){\zeta ^2}\left( 3 \right) - \zeta \left( 2 \right)\zeta _n^2\left( 3 \right)}}{n}}  =  - \frac{{73}}{4}\zeta \left( 2 \right)\zeta \left( 7 \right) + \frac{{21}}{8}\zeta \left( 3 \right)\zeta \left( 6 \right) + 25\zeta \left( 4 \right)\zeta \left( 5 \right) - {\zeta ^3}\left( 3 \right),\\
&\sum\limits_{n = 1}^\infty  {\left( {\bar \zeta \left( 2 \right) - {{\bar \zeta }_n}\left( 2 \right)} \right)\left( {\bar \zeta \left( 3 \right) - {{\bar \zeta }_n}\left( 3 \right)} \right) = } \frac{{21}}{{16}}\zeta \left( 4 \right) + \frac{1}{2}\zeta \left( 2 \right){\ln ^2}2 - \frac{1}{{12}}{\ln ^4}2 - 2{\rm{L}}{{\rm{i}}_4}\left( {\frac{1}{2}} \right) - \frac{3}{8}\zeta \left( 2 \right)\zeta \left( 3 \right),\\
&\sum\limits_{n = 1}^\infty  {\frac{{{\zeta ^2}\left( p \right) - \zeta _n^2\left( p \right)}}{{{n^k}}}}  = \zeta \left( p \right)\zeta \left( {p,k} \right) + \zeta \left( {p,k,p} \right) + \zeta \left( {p,k + p} \right)\quad (k\geq 1,p>1),\\
&\sum\limits_{n = 1}^\infty  {{{\left( {\zeta \left( p \right) - {\zeta _n}\left( p \right)} \right)}^2} = 2} \zeta \left( p \right)\zeta \left( {p,p - 1} \right) + \zeta \left( {2p - 1} \right) - {\zeta ^2}\left( p \right)\quad (p>1).
\end{align*}
These identities can be obtained from the main theorems and corollaries which are presented in the paper.
Some above results are already in the literature, e.g., the fourth, ninth and the last one equations appear as examples of Theorem 4.1 in \cite{MEH2016}.

In fact, by using the method of this paper, it is possible to evaluate other sums involving zeta tails. For example, we have
\begin{align*}
& \sum\limits_{n = 1}^\infty  {\frac{{\zeta \left( {\bar m} \right){\zeta _n}\left( {\bar p} \right) - \zeta \left( {\bar p} \right){\zeta _n}\left( {\bar m} \right)}}
{{{n^k}}}}  = \zeta \left( {\bar p} \right)\zeta \left( {\bar m,k} \right) - \zeta \left( {\bar m} \right)\zeta \left( {\bar p,k} \right), \hfill \\
& \sum\limits_{n = 1}^\infty  {\frac{{\zeta \left( {\bar m} \right){\zeta _n}\left( p \right) - \zeta \left( p \right){\zeta _n}\left( {\bar m} \right)}}
{{{n^k}}}}  = \zeta \left( p \right)\zeta \left( {\bar m,k} \right) - \zeta \left( {\bar m} \right)\zeta \left( {p,k} \right), \hfill \\
&\sum\limits_{n = 1}^\infty  {\frac{{{\zeta ^2}\left( {\bar p} \right) - \zeta _n^2\left( {\bar p} \right)}}
{{{n^k}}}}  = \zeta \left( {\bar p} \right)\zeta \left( {\bar p,k} \right) + \zeta \left( {\bar p,k,\bar p} \right) + \zeta \left( {\bar p,\overline {k + p}} \right),\\
&\sum\limits_{n = 1}^\infty  {\frac{{\zeta \left( {\bar m} \right)\zeta \left( {\bar p} \right) - {\zeta _n}\left( {\bar m} \right){\zeta _n}\left( {\bar p} \right)}}
{{{n^k}}}}  = \zeta \left( {\bar p} \right)\zeta \left( {\bar m,k} \right) + \zeta \left( {\bar p,k,\bar m} \right) + \zeta \left( {\bar p,\overline {m + k}} \right),\\
&\sum\limits_{n = 1}^\infty  {\frac{{\zeta \left( m \right)\zeta \left( {\bar p} \right) - {\zeta _n}\left( m \right){\zeta _n}\left( {\bar p} \right)}}
{{{n^k}}}}  = \zeta \left( {\bar p} \right)\zeta \left( {m,k} \right) + \zeta \left( {\bar p,k,m} \right) + \zeta \left( {\bar p,m + k} \right).
\end{align*}
{\bf Acknowledgments.} The authors would like to thank the anonymous
referee for his/her helpful comments, which improve the presentation
of the paper.

 \end{document}